\documentstyle[12pt]{amsart}
\begin{document}
\newtheorem{Theorem}{Theorem}[section]
\newtheorem{Lemma}[Theorem]{Lemma}
\newtheorem{Corollary}[Theorem]{Corollary}
\newtheorem{Proposition}[Theorem]{Proposition}
\newtheorem{Example}[Theorem]{Example}
\newtheorem{Conjecture}[Theorem]{Conjecture}
\newtheorem*{Definition}{Definition}

\def\core{\operatorname{core}}
\def\height{\operatorname{ht}}
\def\rk{\operatorname{rank}}
\def\Min{\operatorname{Min}}

\def\inc{\subseteq}
\def\To{\longrightarrow}
\def\sk{\smallskip}
\def\mm{{\frak m}}
\def\Im{\text{Im}}
\def\Q{\cal Q}

\title{On the core of ideals}
\author{Craig Huneke}
\date{November 17, 2002}
\address{Department of Mathematics, University of Kansas, Lawrence, KS}
\email{huneke@@math.ukans.edu}
\author{Ng\^o Vi\^et Trung}
\address{Institute of Mathematics,  Box 631, B\`o H\^o, 10000 Hanoi, 
Vietnam}
\email{nvtrung@@math.ac.vn}
\thanks{The first author was partially supported by the NSF. The second 
author is partially supported by the National Basic Research Program of 
Vietnam}
\maketitle

\section*{Introduction} \sk

Let $I$ be an ideal in a Noetherian ring. An ideal $J\subseteq 
I$ is called a {\it reduction} of $I$ if there is a positive number 
$n$ such that $J I^n = I^{n+1}$. In other words, $J$ is a 
reduction of $I$ if and only if $I$ is integrally dependent on 
$J$ [NR]. The {\it core} of $I$, denoted by $\core(I)$, is defined to 
be the intersection of all reductions of $I$. \sk

The core of ideals was first studied by Rees and Sally [RS], partly due 
to its connection to the theorem of Brian{\c c}on and Skoda. Later, 
Huneke and Swanson [HuS] determined the core of integrally closed 
ideals in two-dimensional regular local rings and showed a  close 
relationship to 
Lipman's adjoint ideal. 
Recently, Corso, Polini and Ulrich [CPU1,2] gave 
explicit descriptions for the core of certain ideals in Cohen-Macaulay 
local rings, extending the result of [HuS]. In these two papers, 
several questions and conjectures were raised
which provided motivation  for our work. More recently,
Hyry and Smith [HyS] have shown that the core
and its properties are closely related to a conjecture of Kawamata on 
the
existence of sections for numerically effective line bundles which
are adjoint to an ample line bundle over a complex smooth algebraic
variety, and they generalize the result in [HuS] to arbitrary
dimension and more general rings. 
Nonetheless, there are many unanswered questions on the 
nature of the core. One reason is that it is difficult to 
determine the core and there are relatively few computed examples. \sk

Our focus in this paper is  in effective computation of  the core
with  an eye to partially  answering some questions raised in
[CPU1,2].
A first approach to understanding the core was given by Rees and
Sally. For an ideal $I$ in a local Noetherian ring $(R,\mm)$
having analytic spread $\ell$, one can take $\ell$ generic generators
of $I$ in a ring of the form $R[U_{ij}]_{\mm R[U_{ij}]}$ which
generate an ideal $\Q$. This ideal is then a generic minimal reduction
of the extended ideal and a natural question is whether
$\core(I) = \Q\cap R$. In [CPU1, Thm. 4.7] this equality is proved 
under special conditions,
namely if $R$ is local Cohen-Macaulay ring  with
infinite residue field, and  $I$  
satisfies $G_{\ell}$ and is weakly $(\ell-1)$-residually
$S_2$. While these conditions look somewhat technical, every 
$\mm$-primary
ideal satisfies the conditions. In particular the core of an 
$\mm$-primary ideal
in a local Cohen-Macaulay ring can be computed as the contraction of
a generic reduction. In [CPU1 Example 4.11] they give an example to 
show this is not
true if $R$ is not Cohen-Macaulay.  Under the same conditions, it is
also shown in [CPU1, Thm 4.8] that if $R\rightarrow R'$ is a flat local
homomorphism of Cohen-Macaulay local rings with infinite residue 
fields,
$I$ is an ideal of $R$ of analytic spread $\ell$, $I$ and $IR'$ are
$G_{\ell}$ and universally weakly $(\ell-1)$-residually
$S_2$, then $\core(IR') = \core(I)R'$. They raise as questions whether 
or not
$\core(I)\inc \Q\cap R$ in general and whether or not for arbitrary 
flat
local homomorphisms $\core(I)R'\inc \core(IR')$. We are able to give 
partial answer to both of these questions (see the discussion below).

In [CPU2, Conj. 5.1] a very general conjecture was made concerning how 
to
calculate the core:

 \begin{Conjecture}\label{conjec} Let $R$ be a local Cohen-Macaulay 
ring  with
infinite residue field. Let  $I$  be an $R$-ideal  of analytic
spread  $\ell$ that  satisfies $G_{\ell}$ and is weakly 
$(\ell-1)$-residually
$S_2$. Let $J$ be a minimal reduction of $I$ and let $r$ denote
the reduction  number of $I$ with  respect to $J$.  Then
$$\core(I)= (J^r:I^r)I = (J^r:I^r)J = J^{r+1}:I^r.$$
\end{Conjecture}

An interesting case of this conjecture occurs when $I$ is an 
equimultiple ideal. In  this case the conditions $G_{\ell}$ and weakly 
$(\ell-1)$-residually $S_2$  are automatically  satisfied,
so the  conjecture applies to all such ideals. In  [HyS, Thm. 1.3.3],
it  is shown the core of $I$ is equal  to $J^{n+1}:I^n$
for $n \gg 0$ provided  $I$ is equimultiple, $R$ contains
the rational numbers, and $R[It]$ is Cohen-Macaulay. We are  able to 
verify part of Conjecture~\ref{conjec} for equimultiple ideals in local 
rings with characteristic 0 residue field, namely that $\core(I)= 
J^{r+1}:I^r.$ We have been informed that this result has been obtained 
independently  by Polini and Ulrich  when $\dim R = 1$ or when $R$ is a 
Gorenstein ring [PU]. The one-dimensional case of this conjecture also
follows from the work of Hyry and Smith [HyS].\sk

Another main result of this paper is a closed formula for the graded 
core of the maximal graded ideal of standard graded algebras over a 
field. Let $A$ be a standard graded algebra over a field $k$. 
A {\it reduction} of $A$ is a graded ideal $Q$ generated by linear 
forms such that $Q_n = A_n$ for all large $n$.  Similarly as in the local 
case, we 
define the {\it core} of $A$, denoted by $\core(A)$, to be the 
intersection of the (minimal) reductions of $A$. In other words, 
$\core(A)$ is the intersection of the {\it graded} reductions of the 
irrelevant ideal of $A$.
The importance of the graded core of a homogeneous ideal has recently 
been shown by Hyry and Smith [HyS]. 

If $k$ is an infinite field, every minimal 
reduction of $A$ can be considered as a specialization of an ideal 
$Q_u$ 
generated by generic linear forms in a polynomial ring $A[u]$. 
The parameter space of specializations which are minimal reductions of 
$A$ 
has been described explicitly in [T2]. This leads us to the interesting 
problem of determining the intersection of the specializations of a 
graded module over $k[u]$ upon a given locus. We are able to give an 
effective solution to this problem when $k$ is an algebraically closed 
field (Corollary~\ref{problem}). Combining this result with a 
stratification of the parameter space of minimal reductions of $A$ we 
get a closed formula for $\core(A)$ in terms of the generic ideal $Q_u$ 
(Theorem \ref{core}).
This formula allows us to study basic properties of the graded core. 
For 
instance, we can show that in general, $\core(A \otimes_kE) \neq 
\core(A)\otimes_kE$, where $E$ is a field extension of $k$ (Example 
\ref{question1}). \par

Our formula for the graded core of the maximal irrelevant ideal is 
perhaps more interesting when the base ring is not Cohen-Macaulay. In terms 
of computation, one can always intersect arbitrary minimal reductions 
and hope for stabilization.
This strategy is particularly effective in the Cohen-Macaulay case 
where there is a bound for how many such intersections are needed [CPU1]. 
However, in the non-Cohen-Macaulay case, one has to take into account 
special minimal reductions which can not be chosen arbitrarily, due to 
the stratification in our formula. \sk

The results on the core of graded algebras will be found in Section 1. 
Section 2 largely deals with the core of ideals in local rings. 
We construct counter-examples 
to some open questions on the core raised in [CPU1] and [T2]. In 
particular, we show that the equation $\core(I R') = \core(I)R'$ 
does not hold for an arbitrary flat local homomorphism $R \to R'$ of 
Cohen-Macaulay local rings (Example \ref{question2}).
In Section 3 we settle Conjecture~\ref{conjec} in the affirmative for 
equimultiple ideals in Cohen-Macaulay rings with characteristic zero 
residue field (Theorem \ref{Conjecture}).
In particular, we can prove in the one-dimensional case that  $\core(I) 
= IK$,
where $K$ is the conductor of $R$ in the blowing-up ring at $I$ 
(Theorem \ref{1-dim}). \sk

\noindent{\bf Acknowledgment.} Part of the work of this paper was done 
when the 
second author visited Lawrence, Kansas in the fall of 2001. He would 
like to thank the Department of Mathematics of the University of Kansas 
for support and hospitality. Both authors thank the referee for a careful
reading of this paper.\sk

\section{Core of graded algebras} \sk

Let $A = \oplus_{n\ge 0}A_n$ be a standard graded algebra over an 
infinite field $k$  with $d = \dim A$. Let $x_1,\ldots,x_m$ be linear 
forms which generate the vector space $A_1$. For every point $\alpha = (\alpha_{ij}|\ i = 1,\ldots,d,\ j  =
1,\ldots,m) \in {\Bbb P}_k^N$, $N := dm-1$, we will denote by $Q_\alpha$ the 
ideal generated by $d$ linear forms
$$y_i = \alpha_{i1}x_1 + \cdots + \alpha_{im}x_m\ (i = 1,\ldots,d).$$
We may  consider $Q_\alpha$ as a point of the Grassmannian $G(d,m)$. However, for computational purpose we will study properties of $Q_\alpha$  in terms of  $\alpha$. 
\sk

Let $u = (u_{ij}|\ i = 1,\ldots,d,\ j = 1,\ldots,m)$ be a family 
of $N+1$ indeterminates. Consider $d$ generic elements
$$z_i = u_{i1}x_1 + \cdots + u_{im}x_m\ (i = 1,\ldots,d).$$
For every integer $n \ge 0$ fix a basis ${\cal B}_n$ for the vector 
space $A_n$. For $n \ge 1$ we write every element of the form $z_if$,  
$i = 1,\ldots,d$, $f \in {\cal B}_{n-1}$,  as a linear combination of 
the elements of ${\cal B}_n$. Let $M_n$ denote the matrix of the 
coefficients of these linear  combinations. Then $M_n$ is a matrix with 
entries in $k[u]$. For every integer $t \ge 0$ let $I_t(M_n)$ denote 
the 
ideal of $k[u]$ generated by the $t$-minors of $M_n$. Put $h_n =
\dim_k A_n$. Let $V_n$ be the zero locus of $I_{h_n}(M_n)$ in ${\Bbb 
P}_k^N$.
It was shown in [T2, Theorem 2.1(i)] that ${\Bbb P}_k^N = V_0 
\supseteq V_1 \supseteq \cdots \supseteq V_n \supseteq \cdots$ is a 
non-increasing sequence of projective varieties. This sequence must be 
stationary for $n$ large enough. Hence we can introduce the number
$$r := \max\{n|\ V_n \neq V_{n+1}\}.$$
By [T2, Theorem 2.1(ii)] we get the following parametric 
characterization for the minimal reductions of $A$.

\begin{Proposition} \label{reduction}
$Q_\alpha$ is a minimal reduction of $A$ if and only if $\alpha \not\in 
V_{r+1}$.
\end{Proposition}

The number $r$ introduced above has a clear meaning.
For every minimal reduction $Q$ of $A$ let $r_Q(A)$ denote the largest 
number $n$ such that $Q_n \neq A_n$. Then $r_Q(A)$ is called the {\it =
reduction number}  of  $Q$. If $Q = (x_1,\ldots,x_d)$, then $r_Q(A)$ 
is the maximum degree of the minimal generators of $A$ as a graded 
module over its Noether normalization
$k[x_1,\ldots,x_d]$ [V1]. The supremum of the reductions number of 
minimal reductions of $A$ is called the {\it big reduction number} of 
$A$ [V2] and we will denote it by $br(A)$. It was shown in [T2, 
Corollary 2.3] that
$$br(A) = r.$$
Upper bounds for $br(A)$ in terms of other invariants of $A$ can be 
found in [T1], [V1] and [V2]. 
We note that it is obvious that 
$$\mm^{br(A)+1}\subseteq \core(A).$$

Let $Q_u$ be the ideal of $A[u]$ generated by the generic elements 
$z_1,\ldots,z_d$.
It is clear that every minimal reduction $Q_\alpha$ is obtained from 
$Q_u$ by the substitution $u$ to $\alpha$. If we view $A[u]$ as a 
graded 
algebra over $k[u]$, then every form of a fixed degree $n$ of 
$Q_\alpha$ 
is the evaluation of a form of degree $n$ in $u$ of $Q_u$ at $\alpha$. 
Since $A[u]_n$ is a free module of rank $h_n$ over $k[u]$, we may 
represent every form degree $n$ in $u$ of $Q_u$ as a vector of $h_n$ 
polynomials in $k[u]$.  If we evaluate these polynomials at $\alpha$, 
we 
will get a vector which represents a form of degree $n$ of $Q_\alpha$ 
and all forms of $Q_\alpha$ are obtained in this way. By Proposition 
\ref{reduction},  a form of degree $n$ belongs to $\core(A)$ if and 
only 
if for all $\alpha \not\in V_{r+1}$, it is the evaluation of some form 
of degree $n$ in $u$ of $Q_u$ at $\alpha$.  Therefore, to compute 
$\core(A)$ we need to consider the following problem.\sk

\noindent {\bf Problem}.  Let $E$ be a graded submodule of a free 
module 
$F = k[u]^h$, $h \ge 1$. Consider the elements of $E$ as vectors of 
$h$ polynomials in $k[u]$. Given a projective variety $V \subset {\Bbb P}_k^N$, 
determine the set $$\core_V(E) := \bigcap_{\alpha \not\in 
V}E_\alpha,$$
where $E_\alpha \subseteq k^h$ denotes the vector space generated by 
all 
vectors obtained  of vectors of $E$ by the substitution $u$ to 
$\alpha$. 
\sk

This problem can be effectively solved if $k$ is an algebraically 
closed 
field. In fact our solution can be phrased quite generally and 
computationally solves how to compute $\bigcap_{\mm} \mm M$ where $M$ 
is a finitely
generated $R$-module, $R$ is a Noetherian Jacobson ring, and the 
intersection
runs over all maximal ideals of $R$. Recall that $R$ is said to be 
Jacobson if every prime ideal is the intersection of the maximal ideals 
containing it.
Explicitly:

\begin{Proposition} \label{jacrad} Let $R$ be a Noetherian Jacobson
ring. Let $M$ be a finitely generated $R$-module with a presentation
$$G\overset{\phi}\longrightarrow F\longrightarrow M\longrightarrow 0$$
where $G$ and $F$ are finitely generated free $R$-modules. Let $I_j$ be 
the ideal of $j$-minors of $\phi$, where by convention, $I_0 = R$.
Fix an ideal $J$ in $R$, and let $D(J)\subseteq \text{m-spec}(R)$ be 
the open set in the maximal spectrum of $R$ given by all maximal ideals
not containing $J$.
Then,
$$ \bigcap_{\mm\in D(J)} \mm M = \bigcap_{t\geq 
1}(\sqrt{I_t})M:_MJI_{t-1}.$$
\end{Proposition}

\begin{pf} We first reduce to the case $J = R$.
Let $v\in \bigcap_{\mm\in D(J)} \mm M$. Let $x\in J$ be arbitrary. 
Then $xv\in \bigcap_{\mm\in \text{m-spec}(R)} \mm M$.
We assume we have proved  the proposition in the case that $J = R$, in 
which
case $D(J) = \text{m-spec}(R)$. Then $xv\in \bigcap_{t\geq 
1}(\sqrt{I_t})M:_MI_{t-1}$.
Since $x$ is arbitrary in $J$, it follows that
$v\in\bigcap_{t\geq 1}(\sqrt{I_t})M:_MJI_{t-1}$. Conversely, suppose 
that
$v\in \bigcap_{t\geq 1}(\sqrt{I_t})M:_MJI_{t-1}$. For all $x\in J$, 
$xv\in\bigcap_{t\geq 1}(\sqrt{I_t})M:_MI_{t-1}$, so assuming the case 
in
which $J= R$, we have that $xv\in \bigcap_{\mm\in \text{m-spec}(R)} \mm 
M$.
If $\mm\in D(J)$, then $xv\in \mm M$ will imply that $v\in \mm M$. For,
$M/\mm M$ is a vector space over $R/\mm$, and the image of $x$ in 
$R/\mm$
is a unit. Thus, $v\in  \bigcap_{\mm\in D(J)} \mm M$. It remains to
prove the proposition in the case in which $J = R$, which we henceforth 
assume.

Let $\mm$ be a maximal ideal of $R$. We adopt the following notation:
for a module or element, we write an overline
for the image after tensoring with $R/\mm$. 

We prove that $\bigcap_{t\geq 1}(\sqrt{I_t})M:_MI_{t-1}\subseteq 
\bigcap_{\mm\in \text{m-spec}(R)} \mm M$. Let
$v$ be an  element in $\bigcap_{t\geq 1}(\sqrt{I_t})M:_MI_{t-1}$.
Let $\mm$ be an arbitrary maximal ideal of $R$, and set $r$ equal to 
the rank
of $\overline{\phi}$. Then $I_{r+1}(\phi)\subseteq \mm$, and $I_{r}$ is 
not
contained in $\mm$.
Choose an element $c\notin \mm$ such that $c\in I_{r}$. By assumption,
we have that $cv\in \sqrt{I_{r+1}}M$, and hence 
$\overline{c}\cdot\overline{v} = 0$ in
$\overline{M} = M/\mm M$. But $\overline{c}$ is a nonzero element of 
the field
$R/\mm$, and thus $\overline{v} = 0$, i.e., $v\in \mm M$.

To finish the proof we show the opposite containment, i.e.
$\bigcap_{\mm\in \text{m-spec}(R)}\mm M\subseteq \bigcap_{t\geq 
1}(\sqrt{I_t})M:_MI_{t-1}$.
Let $v\in \bigcap_{\mm\in \text{m-spec}(R)}\mm M$. Choose a lifting 
$u\in F$ of $v$.
By assumption, for all maximal ideals $\mm$, 
$\overline{v} = 0$, so that $\overline{u}$ is in the image of 
$\overline{\phi}$.
We claim that for all $t$, $I_{t-1}v\inc (\sqrt{I_t})M$. To prove this
it suffices to prove that for all $t$,  $I_{t-1}u\inc (\sqrt{I_t})F + 
\Im(\phi)$.
Assume this is not true for some $t$ which we fix. 
Replace $R$ by $R/\sqrt{I_t}$, and $M$ by $M/(\sqrt{I_t})M$.
For all maximal ideals $\mm$ of $R/\sqrt{I_t}$, we have that the image 
of $v$ in
$M/(\sqrt{I_t})M$ is in $\mm (M/(\sqrt{I_t})M)$, and to achieve a 
contradiction,
it suffices to prove that $vI_{t-1} = 0$. Henceforth we assume that
$\sqrt{I_t} = 0$. In this case the rank of $\phi$ is at most $t-1$.
Let $G = R^s$ and $F = R^h$. Adjoin $u$ to the matrix
$\phi$ to get a $h$ by $s+1$ matrix $\psi$ whose last column is $u$. 
For all maximal
ideals $\mm$, the rank of $\overline{\psi}$ is the same as the rank of
$\overline{\phi}$. Hence
the rank of $\overline{\phi}$ and thus the rank of  $\overline{\psi}$ 
is at most $t-1$. From Lemma~\ref{rank} below it follows
that $I_t(\psi) = 0$. Let $e_1,...,e_s$ be the given basis of $G$. If we adjoin  the row $(\phi(e_1),...,\phi(e_s),u)$ to $\psi$ to get a new matrix, then  all new $t$-minors (which involve the row $(\phi(e_1),...,\phi(e_s),u)$) vanish. Expanding these minors we see that the product of any $(t-1)$-minor of $\phi$ with $u$ can be expressed as a linear combination of $\phi(e_1),...,\phi(e_s)$. Therefore, $I_{t-1}(\phi)u \in \Im(\phi)$, proving this direction.
\end{pf}

\begin{Lemma} \label{rank}  Let $R$ be a reduced Jacobson ring, and  
let
$\phi: R^s\longrightarrow R^h$ be a homomorphism of free $R$-modules.
Then $\text{rank}(\phi) = \text{max}\{\text{rank}(\overline{\phi})\}$
where  the maximum is taken over all maximal ideals $\mm$ of $R$,
and where $\overline{\phi}$ denotes  the  map from $(R/\mm)^s
\longrightarrow (R/\mm)^h$  induced by $\phi$.
\end{Lemma}

\begin{pf} Clearly the maximum is at most the rank of $\phi$.
Set $r = \text{rank}(\phi)$. Then the $(r+1)$-size minors
of $\phi$ are zero, and  there is a nonzero $r$ by $r$ minor.
Since $R$ is  reduced and Jacobson, there is a maximal  ideal
$\mm$ which does not contain $I_r(\phi)$. Passing to  $R/\mm$
gives that $r = \text{rank}(\overline{\phi})$.
\end{pf}

Now we use Proposition \ref{jacrad} to study the problem above. Let $E$ 
be a graded submodule of a free module $F = k[u]^h$ for $h\geq 1$. 
Choose a set of generators $g_1,\ldots,g_s$ for $E$. Let $M$ 
denote the $h \times s$ matrix of the coordinates of $g_1,\ldots,g_s$. 
For every integer $t \ge 0$ let $I_t$ be the ideal of $k[u]$ generated 
by the $t$-minors of $M$. We let $M = F/E$, and denote the canonical 
projection of $F$ onto
$M$ by $\pi$.

\begin{Corollary} \label{problem} Assume that $k$ is an algebraically 
closed field. Let $J$ be an ideal of $k[u]$ such that $\sqrt{J}$ is the defining ideal of $V$. Then
$$\core_V(E) = \bigcap_{t \ge 
0}\big[\big(E+\sqrt{I_{t}}F\big):JI_{t-1}\big]\cap k^r.$$
\end{Corollary}

\begin{pf} Let $f$ be an arbitrary vector in $\core_V(E)$. We think of 
$f$ as
living in $k[u]^h$ as a vector of constants.  
Note that $\alpha \notin V$ if and only if the maximal ideal 
$\mm_{\alpha}\in D(J)$.
Then $f\in E_{\alpha}$ if and only if $\pi(f)\in \mm_{\alpha} M$.
Hence by Proposition~\ref{jacrad}, we have that
$f\in \text{core}_V(E)$ if and only if $f\in \bigcap_{\alpha\notin V} 
E_{\alpha}$
if and only if $\pi(f)\in \bigcap_{\mm_{\alpha}\in D(J)} \mm_{\alpha}M 
=
\bigcap_{t\geq 1} (\sqrt{I_t})M:_M JI_{t-1}$ if and only if
$f\in  \bigcap_{t \ge 1}\big[\big(E+\sqrt{I_{t}}F\big):JI_{t-1}\big]$.
\end{pf}

The above corollary allows us to compute the graded pieces of 
$\core(A)$.
But what we need is a closed formula for the whole ideal $\core(A)$. 
Such a formula can be found by taking into account all possible Hilbert 
functions of minimal reductions of $A$. First, we shall describe 
the parameter space of minimal reductions with a given Hilbert 
function.\sk

Let $n$ and $t$ be fixed positive integers. We denote by $V_{n,t}$ the 
zero locus of $I_t(M_n)$ in ${\Bbb P}_k^N$, where $I_t(M_n)$ is the 
ideal of $k[u]$ generated by the $t$-minors of the matrix $M_n$ 
introduced before Proposition \ref{reduction}.
Note that $V_{n,1} \subseteq V_{n,2} \subseteq \cdots V_{n,t}\subseteq 
\cdots$ is a non-decreasing sequence of projective varieties in ${\Bbb 
P}_k^N$ and that $V_n = V_{n,h_n}$, where $h_n = \dim_kA_n$.

\begin{Lemma} \label{Hilbert} With the above notations we have
\begin{align*} & V_{n,t+1} \setminus (V_{n,t}\cup V_{r+1}) =\\
& \{\alpha \in {\Bbb P}_k^N|\ \text{\rm $Q_\alpha$ is a minimal 
reduction of $A$ with $\dim_k(Q_\alpha)_n = t$}\}.
\end{align*}
\end{Lemma}

\begin{pf} Let $M_n(\alpha)$ denote the matrix obtained from $M_n$ by 
the substitution $u$ to $\alpha$.
By the definition of $M_n$ we have $\dim_k(Q_\alpha)_n = \rk 
M_n(\alpha)$. Hence $\dim_k(Q_\alpha)_n = t$ 
since $\alpha \in V_{n,t+1} \setminus V_{n,t}$. By Proposition 
\ref{reduction}, $Q_\alpha$ is a minimal reduction of $A$ if and 
only if $\alpha \not\in V_{r+1}$. The conclusion is 
immediate.
\end{pf}

Let $Q_\alpha$ be an arbitrary minimal reduction of $A$.  For $n \ge 
r+1$ we have $(Q_\alpha)_n = A_n$ by [T2, Theorem 2.1(ii)], hence 
$\dim_k(Q_\alpha)_n$ is independent of the choice of $Q_\alpha$.  
Therefore, the Hilbert function of $Q_\alpha$ is determined by the 
finite sequence of values $\dim_k(Q_\alpha)_n$,  $n = 1,\ldots,r$. \sk

To every sequence $H = \{a_1,\ldots,a_r\}$ of $r$ positive integers we 
associate a set
$$V_H := \bigcap_{n=1}^r \big(V_{n,a_n+1}   \setminus  (V_{n,a_n}\cup 
V_{r+1})\big).$$
We call $H$ an {\it admissible sequence} if $V_H \neq \emptyset$. Let 
$\cal S$ be the set of all admissible sequences. By Lemma \ref{Hilbert} 
we have
\begin{align*}
& V_H = \\
& \{\alpha \in {\Bbb P}_k^N|\ \text{\rm $Q_\alpha$ is a minimal 
reduction of $A$ with $\dim_k(Q_\alpha)_n = a_n$, $n =
1,\ldots,r$}\}.
\end{align*}
Therefore we may view $\cal S$ as the set of all possible Hilbert 
functions of minimal reductions of $\core(A)$. \sk

The next theorem uses the finiteness of $\cal S$ to give a closed 
formula for $\core(A)$.
We believe this is the first general such formula without conditions on 
$A$.

\begin{Theorem} \label{core}
Assume that $k$ is an algebraically closed field. Put $J =
I_{h_{r+1}}(M_{r+1})$. Then
$$\core(A) = \bigcap_{(a_1,\ldots,a_r) \in \cal S}\big[\big(Q_u + 
\sum_{n=1}^r \sqrt{I_{a_n+1}(M_n)}\big):J\prod_{n=1}^rI_{a_n}(M_n) 
\big]\cap A.$$
\end{Theorem}

\begin{pf} For brevity let
$$C := \bigcap_{(a_1,\ldots,a_r) \in \cal S}\big[\big(Q_u + 
\sum_{n=1}^r \sqrt{I_{a_n+1}(M_n)}\big):J\prod_{n=1}^rI_{a_n}(M_n) 
\big]\cap A.$$
We will first show that $A_n \subset C$ for $n \ge r+1$, which 
obviously 
implies that $\core(A)_n \subset C$. Write every element of the form 
$z_if$, $i =
 1,\ldots,d$, $f \in B_r$, as a linear combination of the elements of 
$B_{r+1}$. By definition, $M_{r+1}$ is the matrix of the coefficients 
of 
these linear combinations.
Since $h_{r+1} = \dim_kA_{r+1} = \sharp B_{r+1}$, the product of 
every $h_{r+1}$-minor of
$M_{r+1}$ with an element of $B_{r+1}$ can be written as a linear 
combination of the elements $z_if$. From this it follows that $JA_{r+1} 
\subset Q_u$.
Hence $A_{r+1}  \subset (Q_u:J) \cap A 
\subseteq C.$
Since $A$ is generated by the elements of $A_1$, this implies $A_n 
\subset (A_{r+1}) \subseteq C$ for all $n \ge r+1.$  \par

Now we will show that $\core(A)_n \subset C$ for $n \le r$.  By 
Proposition \ref{reduction} we have
$$\core(A)_n = \bigcap_{\alpha \not\in V_{r+1}}(Q_\alpha)_n.$$
Applying Corollary \ref{problem} we obtain 
$$\bigcap_{\alpha \not\in V_{r+1}}(Q_\alpha)_n \subset \bigcap_{t \ge 
0}\big[\big(Q_u+ \sqrt{I_{t+1}(M_n)}\big):JI_t(M_n)\big]\cap A. $$
For any sequence $(a_1,\ldots,a_r) \in \cal S$ we have
$$\big(Q_u+ \sqrt{I_{a_n+1}(M_n)}\big):JI_{a_n}(M_n)  \subseteq 
\big(Q_u 
+ \sum_{n=1}^r 
\sqrt{I_{a_n+1}(M_n)}\big):J\prod_{n=1}^rI_{a_n}(M_n).$$
Therefore
$$\bigcap_{t \ge 0}\big[\big(Q_u+ 
\sqrt{I_{t+1}(M_n)}\big):JI_t(M_n)\big]\cap A \subseteq C.$$
So we get $\core(A)_n \subset C$ for $n \le r$. Summing up we can 
conclude that $\core(A) \subseteq C$. \par

It remains to show that $\core(A) \supseteq C$. Let $f$ be an arbitrary 
element in $C$. By Proposition \ref{reduction} we have to show that $f
\in Q_\alpha$ for all $\alpha \not\in V_{r+1}$. For every positive 
integer $n \le r$ choose $a_n$ to be the unique positive integer with 
the property $\alpha \in V_{n,a_n+1}\setminus V_{n,a_n}$. Put $H =
\{a_1,\ldots,a_r\}$. Then $\alpha \in V_H$. Hence $H$ is an admissible 
sequence. From this it follows that
$$fJ\prod_{n=1}^rI_{a_n}(M_n) \inc Q_u + \sum_{n=1}^r 
\sqrt{I_{a_n+1}(M_n)}.$$
Note that $\alpha$ is a zero of $\sum_{n=1}^r \sqrt{I_{a_n+1}(M_n)}$ 
and that there exists a polynomial $c(u) \in 
J\prod_{n=1}^rI_{a_n}(M_n)$ such that $c(\alpha) \neq 0$. Since $c(u)f 
\in Q_u+\sum_{n=1}^r \sqrt{I_{a_n+1}(M_n)}$, substituting $u$ to $\alpha$  
we get  $f \in Q_\alpha$. So we have proved that $C \subseteq 
\core(A)$. 
The proof of Theorem \ref{core} is now complete.
\end{pf}

The following example shows that the condition $k$ being  an 
algebraically closed field is necessary in Theorem \ref{core}. 

\begin{Example} \label{residue field} 
{\rm Let $A = {\Bbb R}[x_1,x_2] = 
{\Bbb R}[X_1,X_2]/(X_1^3X_2+ X_1X_2^3,X_2^5)$. Then $\dim A = 1$. Put $Q_u = (u_1x_1+u_2x_2) \subset  {\Bbb R}[u_1,u_2,x_1,x_2]$, where 
$u_1,u_2$ are two indeterminates. For every $n$ let ${\cal B}_n$ be the 
basis of $A_n$ which consists of monomials (in $x_1,x_2$)
which have the possibly highest rank in the lexicographical order. The 
matrices $M_n$ of the coefficients of the elements of the form 
$(u_1x_1+u_2x_2)f$, $f \in {\cal B}_{n-1}$, written as linear 
combinations of elements of ${\cal B}_n$, look as follows:
\begin{align*}
M_1  & = \begin{pmatrix} u_1 & u_2 \end{pmatrix}\\
M_2 & = \begin{pmatrix} u_1 & u_2 & 0\\ 0 & u_1 & u_2 \end{pmatrix}\\
M_3 & = \begin{pmatrix} u_1 & u_2 & 0  & 0 \\ 0 & u_1 & u_2 & 0\\
0 & 0 & u_1 & u_2 \end{pmatrix}\\
M_4 & = \begin{pmatrix} u_1 & u_2 & 0 & 0\\ 0 & u_1 & u_2 & 0\\ 0 & -u_2 & u_1 & 0\\ 0 & -u_1 & 0 & u_2 \end{pmatrix}\\
M_5 & = \begin{pmatrix} u_1 & u_2 & 0 \\  0 & u_1 & u_2\\ 0 & - u_2 & u_1\\ 0 & 0  & -u_1\end{pmatrix}\\
M_6 & = \begin{pmatrix} u_1 & 0 \\  0 & u_2\\ 0 & u_1\end{pmatrix}\\
M_7 & = \begin{pmatrix} u_1 \\  0 \end{pmatrix}\\
M_n & = (u_1)\quad (n 
\ge 8).
\end{align*}
Therefore, $V_1 =V_2 = V_3 = {\Bbb P}^1_ {\Bbb R}$,  $V_4 = \{u_1 = 0\} \cup \{u_2 = 0\}$ and $V_n =
\{u_1=0\}$ for $ n \ge 5$. Hence $r = br(A) = \max\{n|\ V_n \neq 
V_{n+1}\} = 4$ and $J = u_1(u_1,u_2)^2$.  Since the big reduction number is $4$, it
follows that  $(x_1,x_2)^5 \subseteq \core(A)$. Since all minimal reductions of $A$ 
have the form $(\alpha_1x_1+\alpha_2x_2)$ with $\alpha_1 \neq 0$ and 
since $A$ is defined by forms of degree $\ge 4$, we can easily check that $\core(A)$ has no elements of degree $< 4$ and 
$x_1^4,x_1^3x_2,x_1^2x_2^2,x_1x_2^3 \in \core(A)_4$, $x_2^4 \not\in \core(A)_4$. Therefore, $\core(A) = (x_1^4,x_1^3x_2,x_1^2x_2^2)$.
Note that
\begin{align*}
I_t(M_n) & = (u_1,u_2)^t \ \text{if}\ t = 1,...,n\ (n \le 3),\\
I_t(M_4) & = \left\{\begin{matrix}  (u_1,u_2)^t & \text{if} & t = 1,2,3,\\ 
u_1u_2(u_1^2+u_2^2) & \text{if} & t = 4.  \end{matrix}\right.
\end{align*}
Then there are two admissible 
sequence $\{1,2,3,3\}$ and \{1,2,3,4\}.
If Theorem \ref{core} holds for the base field $\Bbb R$, we would 
get
\begin{align*}
\core(A) & = [(Q_u,u_1u_2(u_1^2+u_2^2)):u_1(u_1,u_2)^{11}] \cap 
[Q_u: u_1^2u_2(u_1^2+u_2^2)(u_1,u_2)^8]  \cap A\\
& = (x_1^4 + x_1^2x_2^2,x_1^2x_2^3, x_1x_2^4)  \neq (x_1^4,x_1^3x_2,x_1^2x_2^2),
\end{align*}
which is a contradiction. }
\end{Example}

The next example shows that the formula 
 $\core(A\otimes_kE) = \core(A)\otimes_kE$ does not hold for arbitrary field 
extension $E$ of $k$. 

\begin{Example} \label{question1} {\rm Let $A = {\Bbb R}[x_1,x_2] = 
{\Bbb R}[X_1,X_2]/(X_1^2X_2+X_2^3,X_2^4)$. Then $\dim A = 1$. Put $Q_u = (u_1x_1+u_2x_2) \subset  {\Bbb R}[u_1,u_2,x_1,x_2]$, where 
$u_1,u_2$ are two indeterminates. For every $n$ let ${\cal B}_n$ be the 
basis of $A_n$ which consists of monomials monomials (in $x_1,x_2$)
which have the possibly highest rank in the lexicographical order. The 
matrices $M_n$ of the coefficients of the elements of the form 
$(u_1x_1+u_2x_2)f$, $f \in {\cal B}_{n-1}$, written as linear 
combinations of elements of ${\cal B}_n$, look as follows:
\begin{align*}
M_1  & = \begin{pmatrix} u_1 & u_2 \end{pmatrix}\\
M_2 & = \begin{pmatrix} u_1 & u_2 & 0\\ 0 & u_1 & u_2 \end{pmatrix}\\
M_3 & = \begin{pmatrix} u_1 & u_2 & 0  \\ 0 & u_1 & u_2 \\
0 & -u_2 & u_1 \end{pmatrix}\\
M_4 & = \begin{pmatrix} u_1 & u_2 \\ 0 & u_1 \\ 0 & -u_2 
\end{pmatrix}\\
M_5 & = \begin{pmatrix} u_1 & u_2 \\  0 & u_1 \end{pmatrix}\\
M_6 & = \begin{pmatrix} u_1  \\  0 \end{pmatrix}\\
M_n & = (u_1) \quad (n 
\ge 7).
\end{align*}
So we get $V_1 =V_2 = {\Bbb P}^1_ {\Bbb R}$,  and $V_n =
\{u_1=0\}$ for $ n \ge 3$. Hence $r = br(A) = \max\{n|\ V_n \neq 
V_{n+1}\} = 2$ and $J = (u_1^3+u_1u_2^2)$. As a consequence,
$(x_1,x_2)^3 \subseteq \core(A)$, since in general $\mm^{br(A)+1}\subseteq \core(A)$.
As all minimal reductions of $A$ 
have the form $(\alpha_1x_1+\alpha_2x_2)$ with $\alpha_1 \neq 0$ and 
since $A$ is defined by forms of degree $>2$, one can easily check that 
$\core(A)$ does not contain any form of degree 2. Therefore, $\core(A) 
= (x_1,x_2)^3$.  \par
The core of $A$  will be changed if we replace $\Bbb R$ by $\Bbb C$. In 
this case, we have $V_1 =V_2 = {\Bbb P}^1_ {\Bbb C}$,  $V_3 =  
\{u_1=0\} \cup \{u_1 = \pm iu_2\}$, and $V_n = \{u_1=0\}$ for $ 
n \ge 4$. Hence $r = br(A) =  \max\{n|\ V_n \neq V_{n+1}\} = 3$ and $J = u_1(u_1,u_2)$. Note that 
\begin{align*}
I_t(M_3) & =  (u_1,u_2)^2\ \text{if}\ t =1,2,\\
I_3(M_3) & = u_1(u_1^2+u_2^2).  
\end{align*}
We can easily verify that there are only two admissible sequences 
$(1,2,2)$ and $(1,2,3)$. By Theorem \ref{core} we get
\begin{align*}
\core(A) & = [(Q_u,u_1(u_1^2+u_2^2)):u_1(u_1,u_2)^6] \cap
[Q_u: u_1^2(u_1^2+u_2^2)(u_1,u_2)^4] \cap A \\
&  = (x_1^3 + x_1x_2^2, x_1x_2^3) \neq (x_1,x_2)^3.
\end{align*}}
\end{Example}

The formula of Theorem \ref{core} involves many operations with 
determinantal ideals. Hence we have tried to find a simpler formula. By 
the same argument as in the last part of the proof of Theorem 
\ref{core} 
we always have
$$\big(Q_u:J^\infty) \cap A \subseteq \core(A)$$
where $Q_u:J^\infty$ denotes the set of all elements $f \in A[u]$ such 
that $fJ^n \in Q_u$ for some positive integer $n$.
Since $\sqrt{J}$ is the defining ideal of 
$V_{r+1}$ and ${\Bbb P}_k^N\setminus V_{r+1}$ is the parameter space of 
minimal reductions of $A$, one may raise the question whether $\core(A) 
= 
(Q_u:J^\infty) \cap A$ holds in general. This question has a positive 
answer in the following case.

\begin{Corollary} \label{independence}
Assume that $k$ is an algebraically closed field.  If the Hilbert 
function of every minimal reduction $Q_{\alpha}$ of $A$ does not 
depend on the choice of $\alpha$, then
$$\core(A) = (Q_u:J^\infty) \cap A.$$
\end{Corollary} 

\begin{pf} It is sufficient to show that $\core(A) \subseteq 
(Q_u:J^\infty) \cap A$.
Let $a_n = \rk M_n$, $n = 1,\ldots,r$, and $H =
\{a_1,\ldots,a_r\}$.
Then $I_{a_n}(M_n) \neq 0$ and $I_{a_n+1}(M_n) = 0$. Hence
$V_{n,a_n} \neq {\Bbb P}_k^N$ and $V_{n,a_n+1} = {\Bbb P}_k^N$.
It follows that $$V_H = {\Bbb P}_k^N \setminus 
(\cup_{n=1}^rV_{n,a_n}\cup V_{r+1}) \neq \emptyset.$$
So $H$ is an admissible sequence. By the independence of Hilbert 
functions of minimal reductions, the parameter space of minimal 
reductions equals the parameter space of minimal reductions with 
Hilbert 
functions determined by $H$. Hence we must have $V_H = {\Bbb 
P}_k^N\setminus V_{r+1}$. This implies $V_{n,a_n} \subseteq V_{r+1}$, 
hence $\sqrt{I_{a_n}(M_n)} \supseteq J$ for all $n = 1,\ldots,r$. 
Therefore, applying Theorem \ref{core} we get
$$\core(A) = \big(Q_u:J\prod_{n=1}^rI_{a_n}(M_n)\big) \cap A  
\subseteq  \big(Q_u:J^\infty) \cap A.$$
\end{pf}

The condition on the independence of Hilbert functions of minimal 
reductions is satisfied if $A$ is a Cohen-Macaulay ring. In this case, 
one can even show that
$$\core(A) = (Q_u:J^\infty) \cap A = Q_uA(u) \cap A,$$
where $A(u) := A \otimes_kk(u)$  (the proof is similar to that of 
[CPU1, Theorem 4.7(b)] or [T2, Corollary 4.6]). \sk

Now we will use Theorem \ref{core} to construct a counter-example to 
the 
question whether $\core(A) = (Q_u:J^\infty) \cap A$ holds in 
general.\sk

\begin{Example} \label{question} {\rm Let $A = k[x_1,x_2] 
= k[X_1,X_2]/(X_1^2X_2^2,X_2^5)$. Then $\dim A = 1$. Put $Q_u =
(u_1x_1+u_2x_2) \subset k[u_1,u_2,x_1,x_2]$, where $u_1,u_2$ are two 
indeterminates. For every $n$ let ${\cal B}_n$ be the basis of $A_n$ 
which consists of monomials of degree $n$ in $x_1,x_2$. If we arrange 
these monomials in lexicographical order, then the matrices $M_n$ of 
the 
coefficients of the elements of the form $(u_1x_1+u_2x_2)f$, $f \in 
{\cal B}_{n-1}$, written as linear combinations of elements of ${\cal 
B}_n$, look as follows:
\begin{align*}
M_1  & = \begin{pmatrix} u_1 & u_2 \end{pmatrix}\\
M_2 & = \begin{pmatrix} u_1 & u_2 & 0\\ 0 & u_1 & u_2 \end{pmatrix}\\
M_3 & = \begin{pmatrix} u_1 & u_2 & 0 & 0 \\ 0 & u_1 & u_2 & 0\\
0 & 0 & u_1 & u_2\end{pmatrix}\\
M_4 & = \begin{pmatrix} u_1 & u_2 & 0 & 0\\ 0 & u_1 & 0 & 0\\
0 & 0 & u_2 & 0\\ 0 & 0 & u_1 & u_2 \end{pmatrix}\\
M_5 & = \begin{pmatrix} u_1 & u_2 & 0\\ 0 & u_1 & 0 \\
0 & 0 & u_2\\ 0 & 0 & u_1\end{pmatrix}\\
M_6 & = \begin{pmatrix} u_1 & u_2  \\ 0 & u_1 \\ 0 & 0 \end{pmatrix}\\
M_n & = \begin{pmatrix} u_1 & u_2\\ 0 & u_1 \end{pmatrix}
\quad (n \ge 7).
\end{align*}
So we get $V_1 =V_2 = V_3 = {\Bbb P}^1_k$, $V_4 = \{u_1=0\} \cup 
\{u_2=0\}$, and $V_n = \{u_1=0\}$ for $ n \ge 5$. Hence $r =
br(A) = \max\{n|\ V_n \neq V_{n+1}\} = 4$ and $J = u_1^2(u_1,u_2)$. Note 
that 
\begin{align*}
I_t(M_n) & = (u_1,u_2)^n\ \text{if}\ t = 1,...,n\ (n \le 3),\\
I_t(M_4) & = \left\{\begin{matrix} u_1^2u_2^2 & \text{if} & t = 4,\\ 
(u_1,u_2)^3 & \text{if} & t = 3.  \end{matrix}\right.
\end{align*}
Then we can easily check that there are only two admissible sequences 
$(1,2,3,3)$ and $(1,2,3,4)$. By Theorem \ref{core} we get
\begin{align*}
\core(A) & = [(Q_u,u_1u_2):u_1^2(u_1,u_2)^{10}] \cap
[Q_u: u_1^4u_2^2(u_1,u_2)^7] \cap A\\
& = (x_1^4,x_1^3x_2,x_1x_2^3).
\end{align*}
On the other hand, we have
$$(Q_u:J^\infty) \cap A = (x_1^4,x_1^3x_2,x_1x_2^4) \neq 
(x_1^4,x_1^3x_2,x_1x_2^3).$$}
\end{Example}

One may also ask whether $\core(A) = Q_uA(u) \cap A$ holds for an 
arbitrary graded algebra $A$.  The reason for raising this question is 
the fact that $Q_uA(u)$ is the generic minimal reduction of $A(u)$. As 
noted before, it has a positive answer if $A$ is a Cohen-Macaulay 
ring. But, as in the local case [CPU1, Example 4.11], this question has 
a negative answer in general. \sk

\begin{Example} \label{generic} {\rm
Let $A =k[x_1,x_2]/(x_1x_2,x_2^3)$. Then $\dim A = 1$ and we can 
put $Q_u =(u_1x_1+u_2x_2)$. It is easy to check that $Q_uA(u) \cap A 
= (x_1,x_2)^2$. Since $x_2^2 \not\in (x_1)$, we have $x_2^2 \not\in 
\core(A)$.
Hence $\core(A) \neq Q_uA(u) \cap A$.}
\end{Example}

Following Corso, Polini and Ulrich [CPU1, Question (iii)] one may 
weaken 
the above question  to as whether $\core(A) \subseteq Q_uA(u) \cap A$ 
for an arbitrary graded algebra $A$. Using Theorem \ref{core} we can 
give a positive answer when $k$ is an algebraically closed field.

\begin{Corollary} \label{CPU1}
Assume that $k$ is an algebraically closed field.  Then
$$\core(A) \subseteq Q_uA(u) \cap A.$$
\end{Corollary} 

\begin{pf} Let $a_n = \rk M_n$, $n = 1,\ldots,r$, and $H =
\{a_1,\ldots,a_r\}$. As shown in the proof of Corollary 
\ref{independence}, $H$ is an admissible sequence.
Note that $I_{a_n+1}(M_n) = 0$ for all $n = 1,\ldots,r$. Then we can 
deduce from Theorem \ref{core} that
$$\core(A)A(u) \subseteq \big(Q_u:J\prod_{n=1}^rI_{a_n}(M_n) \big)A(u) 
= Q_uA(u).$$
Since $\core(A) \subseteq \core(A)A(u) \cap A$,  this implies the 
statement.
 \end{pf}

\section{Core of ideals in local rings}

Let $(R,\mm)$ be a local ring with infinite residue field $k$ and $I$ 
an ideal of $R$. Recall that an ideal $J$ is a reduction of $I$ if 
there exists an integer $n$ such that $J I^n = I^{n+1}$. The 
least integer $n$ with this property is called the {\it reduction 
number} of $I$ with respect to $J$, and we will denote it by 
$r_J(I)$. The {\it big reduction number} of $I$, denoted by 
$br(I)$ is the supremum of all reduction numbers $r_J(I)$, where 
$J$ is a minimal reduction of $I$ with respect to inclusion. \sk

These notions are closely related to their counterparts in the graded 
case. Let $F(I) = \oplus_{n\ge 0}I^n/\mm I^n$, the fiber ring of 
$I$. Then $F(I)$ is a standard graded algebra over $k$.  It is known 
that every minimal reduction of $I$ is generated by $d$ elements [NR], 
where  $d = \dim F(I)$ (the analytic spread of $I$). In fact, if 
$c_1,\ldots,c_d$  are elements of $I$ and $c_1^*,\ldots,c_d^*$ are the 
residue classes of $c_1,\ldots,c_d$ in $I/\mm I$, then $J = 
(c_1,\ldots,c_d)$ is a minimal reduction of $I$ if and only if $Q =
(c_1^*,\ldots,c_d^*)$ is a minimal reduction of $F(I)$. Moreover, we 
have $r_J(I) = r_Q(F(I))$ and, in particular, $br(I) =
br(F(I))$  [T2, Lemma 4.1]. \sk

Assume that $I = (a_1,\ldots,a_m)$. Write
$$c_i = \sum_{j=1}^m \alpha_{ij}a_j \  (i =1,\ldots,d)$$
and put $\alpha = (\alpha_{ij}) \in R^{md}$. We may view $J$ as a 
specialization at $\alpha$ of the ideal ${\cal Q} \subset R[u]$ 
generated by the generic elements:
$$b_i = \sum_{j=1}^m u_{ij}a_i\ (i = 1,\ldots,d),$$
where $u = \{u_{ij}|\ i = 1,\ldots,d,\ j = 1,\ldots,m\}$ is a 
family of indeterminates. 
The parameter space of the specializations of ${\cal Q}$ which are 
minimal reductions of $I$ can be described explicitly as follows.\sk

Put $r = br(I)$. Fix two minimal bases of $I^r$ and $I^{r+1}$ 
which consist of monomials in the elements $a_1,\ldots,a_m$. Write the 
elements of the form $b_ig$, where $b_i$ is a generic element of $I$ 
($i = 1,\ldots,d$) and $g$ is a monomial of the fixed basis of 
$I^r$, as a linear combination of the elements of the fixed basis of 
$I^{r+1}$ with coefficients in the polynomial ring $R[u]$. Let $M$ 
denote the matrix of these coefficients.
Let $\cal J$ denote the ideal of $R[u]$ generated by the $h\times 
h$-minors of $M$, where $h$ is the minimal number of generators of 
$I^{r+1}$. It was shown in [T2, Proposition 4.3] that $J$ is a 
minimal reduction of $I$ if and only if ${\cal J}_\alpha = R$, where 
${\cal J}_\alpha$ denotes the ideal of $R$ obtained from $\cal J$ by 
the 
substitution $u$ to $\alpha$. \sk

Unlike the graded case, the above description of the parameter space of 
minimal reductions does not lead to a closed formula for the core. We 
only get the following inclusions.

\begin{Proposition} \label{inclusion}
Let ${\cal Q} :  {\cal J}^\infty$ denote the set of elements $f \in 
R[u]$ with $f{\cal J}^n \subseteq {\cal Q}$ for some positive integer 
$n$. Then 
$$I^{r+1} \subseteq [{\cal Q}:({\cal Q}:I^{r+1})] \cap R \subseteq 
[{\cal Q}:(I^r{\cal Q}:I^{r+1})] \cap R \subseteq ({\cal Q} :  {\cal 
J}^\infty) \cap R \subseteq \core(I).$$
\end{Proposition}

\begin{pf} By the definition of $\cal J$ we have $I^{r+1}{\cal J} 
\subseteq I^r{\cal Q}$. Therefore ${\cal J} \subseteq I^r{\cal 
Q}:I^{r+1}\subseteq {\cal Q}:I^{r+1}$. From this it follows that
$$I^{r+1} \subseteq {\cal Q}:({\cal Q}:I^{r+1}) \subseteq {\cal 
Q}:(I^r{\cal Q}:I^{r+1}) \subseteq {\cal Q}:{\cal J} \subseteq {\cal 
Q}:{\cal J}^\infty.$$
It remains to show that $({\cal Q}:{\cal J}^\infty) \cap R \subseteq 
\core(I)$. But this can be shown by using the above description of the 
parameter space of  minimal reductions of $I$ (see [T2, Theorem 4.4]).
\end{pf}

Corso, Polini and Ulrich [CPU1, Proposition 5.4] have shown that if $R$ 
is a Cohen-Macaulay local ring and $I$ satisfies certain residual 
conditions (which include the case $I$ is an $\mm$-primary ideal), 
then $\core(I) = [{\cal Q} : ({\cal Q}:f^{r+1}) ]\cap R$ for any 
non-zerodivisor $f \in I$. Since ${\cal Q} :({\cal Q}:f^{r+1}) 
\subseteq {\cal Q}:({\cal Q}:I^{r+1})$, this implies
$$\core(I) = [{\cal Q}:({\cal Q}:I^{r+1})]\cap R = [{\cal Q}:(I^r{\cal 
Q}:I^{r+1})] \cap R = ({\cal Q} :  {\cal J}^\infty) \cap R.$$
All the above facts have led  to the question whether $\core(I) = 
({\cal Q} :  {\cal J}^\infty) \cap R$ holds in general [T2, Section 4]. 
Using Example \ref{question} we can now show that this is not the case.

\begin{Example} \label{local-question}
{\rm Let $R = k[x_1,x_2]_{(x_1,x_2)}$ where $k[x_1,x_2]:=
k[X_1,X_2]/(X_1^2X_2^2,X_2^5)$. We will first show that 
$$\core(\mm) = (x_1^4,x_1^3x_2,x_1x_2^3)R.$$ 
Observe that $k[x_1,x_2]$ is the graded algebra $A$ considered 
in Example \ref{question} and that $\core(A) = (x_1^4,x_1^3x_2,x_1x_2^3)$. Since 
every minimal reduction of $A$ generates a minimal reduction of 
$\mm$, we have $\core(\mm) \subseteq \core(A)R = (x_1^4,x_1^3x_2,x_1x_2^3)R$. For the 
converse inclusion, let $(c)$ be an arbitrary minimal reduction of 
$\mm$.
Write $c = x + e$ with $x \in A_1$ and $e \in \mm^2$. Then $x$ 
generates a minimal reduction of $\core(A)$. Hence $(x_1^4,x_1^3x_2,x_1x_2^3) \subseteq xA_3 
= (c-e)A_3$.
Since $r_{(c)}(\mm) \le br(\mm) = br(A) = 4$, we have $eA_3 \subset 
\mm^5 \subseteq (c)$.
Therefore $(c-e)A_3 \subset (c)$, which implies $(x_1^4,x_1^3x_2,x_1x_2^3)R \subset (c)$. So 
we get $(x_1^4,x_1^3x_2,x_1x_2^3)R \subseteq \core(\mm)$. On the other hand, if we put 
${\cal Q} = (u_1x_1+u_2x_2)R$, where $u_1,u_2$ are two indeterminates, 
then the matrix $M$ is the matrix $M_5$ of Example \ref{question}:
$$M  = \begin{pmatrix} u_1 & u_2 & 0\\ 0 & u_1 & 0 \\
0 & 0 & u_2\\ 0 & 0 & u_1\end{pmatrix}$$
Since $\mm^5$ is generated by 3 elements, we get ${\cal J} = I_3(M) 
= u_1^2(u_1,u_2)$. Now it can be easily checked that
$$({\cal Q} :  {\cal J}^\infty) \cap R = (x_1^4,x_1^3x_2,x_1x_2^4)R 
\neq (x_1^4,x_1^3x_2,x_1x_2^3)R = \core(\mm).$$}
\end{Example}

It is a natural problem to extend the results on the graded core to the 
core of maximal ideals in local rings which arise from graded algebras. 
So we raise the following basic question:\sk

\noindent{\bf Problem.} Let $(R,\mm)$ be the localization of a standard 
graded algebra $A$ over an infinite field at its maximal graded ideal. 
Is it 
true that $\core(\mm) = \core(A)R$?\sk

Note that the above formula holds for the local ring of Example 
\ref{local-question}. Now we shall see that this formula also holds if 
$br(A) \le 2$.

\begin{Proposition} Let $(R,\mm)$ be the localization of a standard 
graded algebra $A$ over an infinite field at its maximal graded ideal. 
Assume that $br(A) \le 2$. Then
$$\core(\mm) = \core(A)R.$$
\end{Proposition}

\begin{pf} Since every minimal reduction of $A$ generates a 
minimal reduction of $\mm$, we always have $\core(\mm) \subseteq 
\core(A)R$. For the converse inclusion, let $J = (c_1,\ldots,c_d)$ 
be an arbitrary minimal reduction of $\mm$, where $d := \dim A$.
Write $c_i = x_i + e_i$ with $x_i \in A_1$ and $e_i \in \mm^2$. Let 
$Q$ be the ideal generated by the elements $x_1,\ldots,x_d$. Then $Q$ 
is 
a minimal reduction of $A$. Hence
$$\core(A)_2 \subseteq Q_2 \subseteq J + (e_1,\ldots,e_d)\mm \subseteq 
J+\mm^3.$$
Since $br(\mm) = br(A) = 2$, we have $\mm^3 \subseteq J$.
Therefore $\core(A)_2 \subset J$. This implies $\core(A) \subset J$ 
because $\core(A)$ is generated by elements of degree $\ge 2$ and 
$\core(A)_3 \subset \mm^3 \subseteq J$. So we can conclude that 
$\core(A)R \subseteq \core(\mm)$.
\end{pf}

Corso, Polini and Ulrich have asked whether $\core(I R') \supseteq 
\core(I)R'$ does hold for every flat local homomorphism $R \to R'$ of 
Cohen-Macaulay local rings.
We shall use the results in the graded case to construct a 
counter-example to this question.

\begin{Example} \label{question2}
{\rm Let $R = {\Bbb R}[x_1,x_2,x_3,x_4]_{(x_1,x_2,x_3,x_4)}$, where
$${\Bbb R}[x_1,x_2,x_3,x_4]= {\Bbb 
R}[X_1,X_2,X_3,X_4]/(X_1^2X_2+X_2^3-X_1^3X_3,X_2^4-X_1^4X_4).$$ 
Then $R$ is a 
two-dimensional Cohen-Macaulay local ring.
Let $I = (x_1,x_2)$. Since $R$ is defined by equations which are 
homogeneous in $x_1,x_2$, we can easily check that
$$F(I) \cong {\Bbb R}[X_1,X_2]/(X_1^2X_2+X_2^3,X_2^4),$$
which is the graded algebra of Example \ref{question1}. From this it 
follows that $br(I) = br(F(I)) = 2$. Hence $I^3 \subseteq 
\core(I)$. Moreover, every minimal reduction of $I$ is generated by 
an element and this element has precisely the form $x = \alpha_1x_1 + 
\alpha_2x_2 + y$ with $\alpha_1 \neq 0$ and $y \in \mm I$. Thus, every 
element of $\core(I)$ can be written both as $ax_1$ and $b(x_1+x_2)$ 
for some elements $a, b \in R$, and we will obtain a relation of the 
form $ax_1-b(x_1+x_2) = 0$. Since $R$ is defined by homogeneous 
equations  of degree $ \ge 3$ in $x_1,x_2$, we can easily verify that 
$\core(I)$ is contained in $I^3$. So we get $\core(I) =
I^3$.\par
Let $R' = R \otimes_{\Bbb R}{\Bbb C}$.  Similarly as above, we can 
also show that $$F(I R') \cong {\Bbb 
C}[X_1,X_2]/(X_1^2X_2+X_2^3,X_2^4)$$
and that $\core(I R')$ is contained in $I^3 R'$. Assume that 
$\core(I R') = I^3R'$. For degree reason, we must have
$I^3 R' = JI^2$ for any minimal reduction $J$ of $I R'$.
By the relationship of the minimal reductions of $I R'$ and $F(I 
R')$, this implies  $ \core(F(I R')) = (X_1,X_2)^3F(I R')$. 
But we have seen in Example \ref{question1} that 
$$\core(F(I R')) \neq (X_1,X_2)^3F(I R').$$
So we get a contradiction. Now we can conclude that $\core(I R')$ is 
strictly contained in $\core(I)R' = I^3R'$.}
\end{Example}

\section{Core of equimultiple ideals in Cohen-Macaulay local rings}

In this section we will concentrate on the core of equimultiple ideals 
in 
 Cohen-Macaulay local rings with characteristic zero residue field.
Recall that an ideal $I$ is {\it equimultiple} if the minimal 
reductions of $I$ are generated by $h$ elements, where $h = \height(I)$. Our aim 
is to prove Conjecture~\ref{conjec} for this case. \par

First, we will consider the one-dimensional case. We shall need the 
following observation on the ubiquity of minimal reductions.

\begin{Lemma} \label{line}
Let $(R,\mm)$ be an one-dimensional local ring with infinite residue field. Let $I$ be an $\mm$-primary ideal. Let $x$ be an 
element of $I$ which generates a minimal reduction of $I$ and $y$ an 
arbitrary element in $I \setminus \mm I$. Then we can find infinitely many 
units $u \in R$ with different residue classes in $R/\mm$ such that 
$(x-uy)$ is a minimal reduction of $I$.
\end{Lemma}

\begin{pf} 
Every element $z \in I \setminus \mm I$ can be 
represented by a point $\alpha$ in a projective space over the field $k = 
R/\mm$. Choose a set of generators $x_1,\ldots,x_m$ of $I$ and write $z = 
u_1x_1+ \cdots + u_mx_m$, where $u_1,\ldots,u_m$ are elements of $R$. 
Let $\alpha_i$ denote the image of $u_i$ in the residue field $k = 
R/\mm$. By Proposition \ref{reduction}, there exists a proper projective 
variety $V \subset {\Bbb P}_k^{m-1}$ such that $(z)$ is a minimal reduction 
of $I$ if and only if $(\alpha_1,\ldots,\alpha_m) \not\in V$.  If $x,y$ 
are represented by the points $\alpha, \beta \in {\Bbb P}_k^{m-1}$, 
then $x-uy$ is represented by the point $\alpha- c\beta$, where $c$ 
denotes the image of $u$ in $k$. The conclusion follows from the fact that 
the line $\alpha -c\beta$ intersects $V$  only at finitely many points.
\end{pf}

Using the above lemma we are able to characterize the core of 
$\mm$-primary ideals in an one-dimensional Cohen-Macaulay local ring $R$ as 
follows. For a given ideal $I$ of $R$ we will denote by $B(I)$ the 
blowup of  $I$, that is, $B(I) = \cup_{n \ge 
1}(I^n:_FI^n)$, where $F$ denotes the ring of fractions of $R$.
The ring $B(I)$ is always in the integral closure of $R$.

\begin{Theorem}  \label{1-dim}
Let $(R,\mm)$ be a one-dimensional Cohen-Macaulay local ring whose 
residue field has characteristic $0$. Let $I$ be an $\mm$-primary ideal, 
and 
let $B$ be the blowup of $I$. Set $K$ equal to the conductor of
$B$ into $R$, and let $x$ be an arbitrary minimal reduction of $I$.
Then $\core(I) = xK = IK$.
\end{Theorem}

\begin{pf}
By [BP, Theorem 1], we have that $IK = zK$ for every element $z \in R$ 
which
generates a minimal reduction of $I$. This implies that $xK = IK 
\subseteq \core(I)$.
Since $K = \cap_{r \ge 1}(x^r):I^r$, it remains to prove that for any 
element $ax \in \core(I)$ and any $r \ge 1$, we have $a \in (x^r):I^r$. 
Since the residue field has characteristic $0$, $I^r$ is generated by 
the $r$-th powers of the elements of $I$ (see. e.g. [Hu, Exercise 
1.11]). Therefore,  it is sufficient to show that $ay^r \in (x^r)$ for any 
element $y \in I \setminus \mm I$. By Lemma \ref{line}, there exist 
infinitely many units $u \in R$ such that with different residue classes in 
$R/\mm$ such that $(x-uy)$ is a minimal reduction of $I$. For such $u$ 
we can write
$$ax = b_u(x-uy)$$
for some element $b_u \in R$. Set $c = y/x$. Then $c \in B$ [No, Lemma 
1]. Hence, for $r$ large enough, there is an integral relation 
$$d_0 + d_1c + \cdots + d_{r-1}c^{r-1} + c^r = 0$$
with $d_0,...,d_{r-1} \in R$. Since $a = b_u(1-uc)$, we have
$$a(1+uc+ \cdots + u^{i-1}c^{i-1}) = b_u(1-u^ic^i)$$
for all $i \ge 1$. It follows that
\begin{align*}
& \sum_{i=0}^{r-1}d_iu^{r-i}[a(1+uc+ \cdots + u^{i-1}c^{i-1})\big] + 
a(1+uc+ \cdots + u^{r-1}c^{r-1}) \\
& = \sum_{i=0}^{r-1}d_iu^{r-i}b_u(1-u^ic^i) + b_u(1-u^rc^r) = 
\sum_{i=0}^{r-1}d_iu^{r-i}b_u + b_u \in R.
\end{align*}
So we get the relation
\begin{align*}
(d_0u^r + d_1u^{r-1} + \cdots + d_{r-1}u + 1)a &+  (d_1u^{r-1} + \cdots 
+ d_{r-1}u + 1)uac +\\  \cdots  & + (d_{r-1}u + 1)u^{r-2}ac^{r-2} + 
u^{r-1}ac^{r-1} \in R.
\end{align*}\par
Varying $u$ we obtain linear equations in
$a,ac,..., ac^{r-2},ac^{r-1}$ modulo $R$. The rows of the coefficient 
matrix of $r$ such equations have the form
$$(d_0u^r + d_1u^{r-1} + \cdots + d_{r-1}u + 1), (d_1u^{r-1} + \cdots + 
d_{r-1}u + 1)u,...,(d_{r-1}u + 1)u^{r-2},u^{r-1}.$$
Decomposing the $(r-1)$th column we see that the determinant of this 
matrix is equal to that of the matrix of the rows
$$(d_0u^r + d_1u^{r-1} + \cdots + d_{r-1}u + 1), (d_1u^{r-1} + \cdots + 
d_{r-1}u + 1)u,...,u^{r-2},u^{r-1}.$$
Decomposing the other columns we will come to the Vandermonde matrix of 
the rows $1,u,...,u^{r-2},u^{r-1}$.
By Lemma \ref{line}, we may choose the elements $u$ in such a way that
the determinant of the Vandermonde matrix is a unit in $R$, and this
proves that $ac, ac^2,...,ac^{r-1} \in R$. Since $r$ can be any number 
large enough, we get $ac^r\in R$ for all $r \ge 1$. This implies $ay^r 
\in (x^r)$, as desired.
\end{pf}

The above theorem allows us to give a positive answer to 
Conjecture~\ref{conjec} in the special case in which $R$ is as in the above theorem. 
This is due to the following set of equivalences:

\begin{Proposition}\label{equivalences} Let $(R,\mm)$ be a 
one-dimensional local ring
with infinite residue field, and let $I$ be an ideal containing a 
nonzerodivisor.
Let $B$ be the blowup of $I$ and let $K$ be the conductor of $B$ into 
$R$.
Let $x$ be an arbitrary minimal reduction of $I$ with reduction number 
$r$.
Then 
$$xK = ((x^r):I^r)I  = ((x^r):I^r)x = ((x^{r+1}):I^r),$$ 
and so the following are equivalent:
\newline\noindent {\rm (i)} $\core(I) = xK$.
\newline\noindent {\rm (ii)} $\core(I) = ((x^r):I^r)I$.
\newline\noindent {\rm (iii)} $\core(I) = ((x^r):I^r)x$.
\newline\noindent {\rm (iv)} $\core(I) = ((x^{r+1}):I^r)$.
\end{Proposition}

\begin{pf} Observe that every minimal reduction of $I$, in particular 
$x$, is a 
nonzerodvisor. We first prove that $((x^r):I^r)I  = ((x^r):I^r)x = 
((x^{r+1}):I^r)$.
Clearly $((x^r):I^r)x\inc ((x^r):I^r)I$ and $((x^r):I^r)x\inc 
((x^{r+1}):I^r)$.
We prove that $((x^{r+1}):I^r)\inc ((x^r):I^r)I\inc ((x^r):I^r)x$ to 
prove all
three are equal. Let $s\in ((x^{r+1}):I^r)$. Then in particular, $sx^r 
= tx^{r+1}$
for some $t\in R$, so that as $x$ is a nonzerodivisor, $s = xt$. As
$xtI^r\inc (x^{r+1})$, we have that $t\in ((x^r):I^r)$, showing that 
$s\in
((x^r):I^r)x\inc ((x^r):I^r)I$. Finally, let $u\in ((x^r):I^r)I$. Then 
$x^{r-1}u\in (x^r)$
which forces $u = xv$ for some $v\in R$. But then $xvI^r\inc I^r\cdot 
I((x^r):I^r) =
xI^r((x^r):I^r)$ (because $I^{r+1} = xI^r$), so that $v\in 
((x^r):I^r)$, and
$u\in ((x^r):I^r)x$. 

It remains to prove that $xK$ is equal to any of the other three 
ideals. We prove that
$K = ((x^r):I^r)$. If $b\in K$,
then $bI^n\inc (x^n)$ for all $n\geq 1$, since $I^n/x^n\inc B$. 
Conversely, if $w\in ((x^r):I^r)$, then we claim that
$w\in K$. For let $b\in B$. We can write $b = a/x^n$ for some $a\in 
I^n$.
Without loss of generality we can assume that $n\geq r+1$. Then $wa\in 
wI^n\inc
wx^{n-r}I^{r}\inc (x^n)$, which implies that $wb\in R$. Thus $K = 
((x^r):I^r)$,
finishing the proof.
\end{pf}

\begin{Corollary}\label{cpu-conj}  Let $(R,\mm)$ be a one-dimensional
Cohen-Macaulay local ring whose residue field has characteristic $0$. 
Let $I$ be an $\mm$-primary ideal of $R$, let $x$ be a minimal
reduction of $I$, and let $r$ be the reduction number of $I$ with 
respect to $x$.
Then 
$$\core(I) = ((x^r):I^r)I  = (x^r:I^r)x = ((x^{r+1}):I^r).$$
\end{Corollary}

\begin{pf} The corollary follows at once from Theorem~\ref{1-dim} and 
Proposition \ref{equivalences}.
 \end{pf}

It is easy to find  examples which show that Corollary \ref{cpu-conj} 
does not hold if the ring $A$ is not Cohen-Macaulay.

\begin{Example}
{\rm Let $A = k[x_1,x_2]_{(x_1,x_2)}$ with $k[x_1,x_2] := 
k[X_1,X_2]/(X_1X_2,X_2^2)$. Let $I$ be the maximal ideal of $A$.
Then $(x_1)$ is a minimal reduction of $I$ with reduction number  $1$. 
We have $(x_1^2):I = (x_1,x_2).$ Since  $x_2$ is not contained in 
$(x_1)$, we must have $\core(I) \neq (x_1^2):I$.}
\end{Example}

To pass from the one-dimensional case to the general case we shall need 
the following result of Hyry and Smith.

\begin{Lemma} \label{HyS} {\rm [HyS, Lemma 5.1.3]}
Let $R$ be a Cohen-Macaulay local ring. Let $I$ be an equimultiple 
ideal of $R$ with $h = \height(I) \ge 1$, let $J$ be a minimal
reduction of $I$, and let $r$ be the reduction number of $I$.
Then the ideal $J^{r+1}:I^r$ does not depend on the choice of $J$.
\end{Lemma}

Now we are able to settle Conjecture \ref{conjec} for equimultiple 
ideals in Cohen-Macaulay local rings whose residue field has characteristic 
$0$.

\begin{Theorem} \label{Conjecture}
 Let $R$ be a Cohen-Macaulay local ring whose residue field has 
characteristic $0$. Let $I$ be an equimultiple ideal of $R$ with $h = 
\height(I) \ge 1$, let $J$ be a minimal
reduction of $I$, and let $r$ be the reduction number of $I$.
Then 
$$\core(I) = J^{r+1}:I^r.$$
\end{Theorem}

\begin{pf} 
Since $J$ is generated by a regular sequence, we have
$J^{r+1}:I^r \subseteq J^{r+1}:J^r = J.$
Since $J^{r+1}:I^r$ does not depend on $J$ by Lemma \ref{HyS}, this 
implies
$J^{r+1}:I^r \subseteq \core(I).$
It remains to show that $\core(I) \subseteq J^{r+1}:I^r$. This 
inclusion is equivalent to the formula  
$$\core(I) \subseteq J^{n+1}:I^n$$ for all $n \ge 0$. In fact, we 
always have 
$J^{r+1}:I^r \subseteq J^{r+1}:J^{r-n}I^n = J^{n+1}: I^n$
for $n \le r$ and  $J^{n+1}:I^n = J^{n+1} : J^{n-r}I^r = J^{r+1}:I^r$ 
for $n > r$. We will use induction on $n$ and $h$ to prove the above 
formula. \par

If $n = 0$, it is trivial because $\core(I) \subseteq J = J:I^0$.
If $h = 1$, we will go back to the one-dimensional case. Let 
$\Min(J^{r+1})$ denote the set of the minimal associated primes of $J^{r+1}$. 
Since 
$J^{r+1}$ is an unmixed ideal, we have
$$J^{r+1} = \bigcap_{P \in \Min(J^{r+1})}J^{r+1}A_P \cap A.$$
By [CPU1, Theorem 4.5] we have
$$\core(I) \subseteq \bigcap_{P \in \Min(J^{r+1})}\core(I_P) \cap A.$$
Since $\dim A_P = 1$, we may apply Corollary \ref{cpu-conj} to 
$\core(I_P)$ and obtain
\begin{align*}
\core(I) & \subseteq  \bigcap_{P \in \Min(J^{r+1})}(J^{r+1}A_P:I^r) 
\cap A\\
& \subseteq \big(\bigcap_{P \in \Min(J^{r+1})}J^{r+1}A_P \cap 
A\big):I^r = J^{r+1}:I^r.
\end{align*}
As noted above, this implies the formula $\core(I) \subseteq 
J^{n+1}:I^n$  for all $n \ge 0$ when $h = 1$.\par

Now let $n > 0$ and $h > 1$. Using induction on $n$ we may assume that 
$\core(I) \subseteq J^n:I^{n-1}$.
Let $S$ denote the set of all superficial elements of $I$ which belong to  
minimal bases of $J$. Note that all elements of $S$ are necessarily non-zerodivisors.
For an arbitrary $x \in S$ let $\bar R = R/(x)$, 
$\bar I = I/(x)$ and $\bar J = J/(x)$. Then $\bar A$ is a 
Cohen-Macaulay local ring, $\bar I$ is an equimultiple ideal with $\height \bar I = 
h-1$, and $\bar J$ is a minimal reduction of $\bar I$. Using induction 
on $h$ we may assume that
$\core(\bar I) \subseteq \bar J^{n+1}:\bar I^n$. Let $a_1,..., a_{h-1}$ be elements of $I$ such that $(a_1,...,a_{h-1},x)/(x)$ is a minimal reduction of $\bar I$. Then
$I^{n+1} \subseteq (a_1,...,a_{h-1})I^n + (x)$
for $n \gg 0$.  Since $x$ is a superficial element of $I$, there exists an integer $c$ such that $(I^{n+1}:x) \cap I^c = I^n$ for $n \gg 0$. By Artin-Rees lemma, $(x) \cap I^{n+1} \subseteq xI^{n-n_0}$ for some integer $n_0$. Hence
$$(x) \cap I^{n+1}= x[(I^{n+1}:x) \cap I^{n-n_0}] = x[(I^{n+1}:x) \cap I^c \cap I^{n-n_0}] = xI^{n}$$
for $n \gg 0$.  These facts imply
$$I^{n+1} = (a_1,...,a_{h-1})I^n + (x) \cap I^{n+1} = (a_1,...,a_{h-1},x)I^n$$
for $n \gg 0$. So $(a_1,...,a_{h-1},x)$ is a minimal reduction of $I$. Now it is clear that  $\core(\bar I) = \cap K/(x)$, where $K$ runs through 
all minimal reductions of $I$ containing $x$. Therefore,
$\cap K \subseteq (J^{n+1},x):(I^n,x)$. Since $\core(I) \subseteq \cap 
K$, we get
$$\core(I)I^n \subseteq J^n \cap (J^{n+1},x) = xJ^{n-1} + J^{n+1}.$$
Let $G = \oplus_{t \ge 0}J^t/J^{t+1}$ be the associated graded ring of 
$J$.
Then $G$ is a Cohen-Macaulay ring with  $\height G_+ \ge 2$. 
Since $S$ corresponds to the complement of a finite union of linear subpaces of $J/\mm J$, we can  choose elements $x_1,...,x_{n+1}  \in S$ such that their initial forms $x_1^*,...,x_{n+1}^*$ in $G$ are pairwise relatively prime. We 
have  $(x_1^*)\cap ... \cap (x_{n+1}^*) = (x_1^*...x_{n+1}^*)$. Since 
$(x_1^*...x_{n+1}^*)_n = 0$, we get
$$\core(I) \subseteq  (x_1J^{n-1}+J^{n+1}) \cap ... \cap 
(x_{n+1}J^{n-1}+J^{n+1}) = J^{n+1}.$$
This proves the formula $\core(I) \subseteq J^{n+1}:I^n$ for all $n > 
0$ and $h > 1$.
\end{pf}

\section*{References}

\noindent [BP] V.~Barucci and K. Pettersson, On the biggest maximally 
generated ideal as the conductor in the blowing up rings, Manuscripta 
Math. 88 (1995), 457-466. \par

\noindent [CPU1] A.~Corso, C.~Polini and B. Ulrich, The structure of 
the 
core of ideals, Math. Ann. 321 (2001), 89-105 \par

\noindent[CPU2] A.~Corso, C.~Polini and B. Ulrich, Core and residual 
intersections of ideals, Trans. Amer. Math. Soc.  354 (2002), 
2579-2594. \par

\noindent [Hu] C.~Huneke, Tight closures and its applications, CBMS 
Regional Conference Series in Math. 88, Amer. Math. Soc., Providence, 
1996.\par

\noindent [HuS] C. Huneke and I. Swanson, Cores of ideals in 
2-dimensional regular local rings, Michigan Math. J. 42 (1995), 
193-208. 
\par

\noindent[HyS] E. Hyry and K.E. Smith, 
On a non-vanishing conjecture of Kawamata and the core of an ideal,
preprint, 2002. \par

\noindent [Na] M. Nagata, Local rings, Interscience, New York, 
1960.\par

\noindent [No] D. G. Northcott, On the notion of a first neighbourhood 
ring, Proc. Cambridge Philos. Soc. 53 (1959), 43-56.\par

\noindent [NR] D. G. Northcott and D. Rees, Reductions of ideals
in local rings, Proc. Cambridge Philos. Soc. 50 (1954), 145-158.\par

\noindent [PU] C. Polini and B. Ulrich, A conjecture on the core of 
ideals, work in progress. \par

\noindent [RS] D. Rees and J. Sally, General elements and joint 
reductions,
Michigan Math. J. 35 (1988), 241-254. \par

\noindent [T1] N.~V.~Trung, Reduction exponent and degree bound
for the defining equations of graded rings, Proc. Amer. Math.
Soc. 101 (1987), 229-236. \par

\noindent [T2] N.~V.~Trung, Constructive characterizations of the 
reduction numbers, Compositio Math. 137 (2003), 99-113. \par

\noindent [V1] W. Vasconcelos, The reduction number of an algebra,
Compositio Math. 106 (1996), 189-197. \par

\noindent [V2] W. Vasconcelos, Reduction numbers of ideals,
J. Algebra 216 (1999), 652-664. \par
\end{document}